\documentclass[12pt]{amsart}
 
\theoremstyle{definition}

\theoremstyle{parrafo}

\numberwithin{equation}{theorem}
 
\begin{document}

\title[]{An example on
the maximal function associated to a nondoubling measure}
 
\author{J. M. Aldaz}
\address{Departamento de Matem\'aticas y Computaci\'on,
Universidad  de La Rioja, 26004 Logro\~no, La Rioja, Spain.}
\email{aldaz\@dmc.unirioja.es}


\thanks{Partially supported by Grant BFM2003-06335-C03-03 of the D.G.I. of
Spain}
 
 
 
 
 

\begin{abstract} We show that there is a measure $\mu$, defined  on the
hyperbolic plane and with polynomial  growth, such that the centered maximal
operator associated to $\mu$ does not satisfy weak type $(1,1)$ bounds.
\end{abstract}

 
\maketitle
 
 
\section
{Introduction and main result}

\markboth{J. M. Aldaz}{The centered maximal operator}

Let $X$ be a metric space and let $\mu$ be a Borel measure defined on $X$. In a
paper (\cite{NTV}) that has exerted considerable influence in later
developments, F. Nazarov, S. Treil, and A. Volberg
showed that a good deal of the theory of Calder\'on Zygmund operators still
holds if one replaces the doubling condition on the measure $\mu$ by the
following polynomial  growth condition: There exist constants $c,d$ such that
for every 
$x\in X$ and every radius $r>0$, $\mu (B(x, r)) \le c r^d$. 
Polynomial growth is a natural assumption in this
area: In the euclidean
case $X= \Bbb R^{d}$, G. David
showed  it is needed for
the $L^2$ boundedness of singular integral operators that, like the
Hilbert transform, are associated to kernels $K$ satisfying  
$|K(x,y) |\ge C|x-y|^{-d}$ (cf. Comment 1, page 60 of \cite{Pa}). 

The ``Calder\'on Zygmund philosophy" consists in controlling
singular integral operators via the maximal function, which
for doubling measures is of weak type $(1,1)$. 
In this context it is natural to ask whether the assumption
$\mu (B(x, r)) \le c r^d$ can replace the doubling hypothesis
and still yield weak type $(1,1)$ bounds. F. Nazarov, S. Treil, and A. Volberg
bypassed this issue in (\cite{NTV}), resorting instead to a
modified maximal function $\tilde M$, where one does not take the usual average
but
divides by the measure of the ball with the same center and
triple radius: $\tilde
Mf(x):=\sup_{r>0}(\mu(B(x,3r)))^{-1}\int_{B(x,r)}|f|\,d\mu$. This modified
maximal function is of weak type $(1,1)$, but in principle does not control
anything.  Polynomial growth is then used to show that there exist ``large
doubling balls" centered at every point, making it possible to utilze 
$\tilde M$ in order to find bounds for other operators.

 A positive answer to the question whether
the hypothesis
$\mu (B(x, r)) \le c r^d$ is sufficient to ensure that
the  weak type $(1,1)$ of the centered maximal function 
would have allowed a development of the theory more
in parallel with the classical case. The negative
answer we obtain shows that such reduction 
is not possible. 
In the example we present the metric space $X$ is just
the hyperbolic plane $\Bbb H$, with a suitably defined 
Borel measure.
More precisely, we prove the following 

\   

\noindent{\bf Theorem.} {\em There exists a Borel measure
$\mu$ on the hyperbolic plane $\Bbb H$ and a constant $c > 0$ 
such that for every $w\in \Bbb H$ and every radius $s>0$, $\mu (B(w, s)) \le c
s$, and the centered maximal function associated to $\mu$ is not of
weak type $(1,1)$.}

\
 
Let us make a comment on terminology: In \cite{NTV} and some later
papers,  a measure satisfying
$\mu (B(x, r)) \le c r^d$ is called $d$-dimensional. But 
the preceding condition does not really give a genuine notion of dimension:
some measures can have uncountably many ``dimensions" (for instance,
planar Lebesgue measure on the unit disc), while other perfectly
good measures have none (for example, the doubling measure 
$\mu (A) := \int_A
x^2 dx$ on $\Bbb R$,  which under any
reasonable definition ought to be regarded as one-dimensional).
So we think it is more precise to
speak of polynomial growth, as is done, for instance, in \cite{Pa}.

 This research was carried out while visiting the Universidad Aut\'onoma de
Madrid. I am indebted to the Department of Mathematics,
 and specially to Prof. Jos\'e Garc\'{\i}a-Cuerva, for the invitation,
 and also for several
useful conversations.

\section{Proof of the theorem}
 
We shall utilize the upper half plane model of 
 the hyperbolic plane $\Bbb H$. The following properties of this model,
which can be found in geometry textbooks (see for instance,
\cite{McC} page 237), will be used in the proof.
The hyperbolic balls $B=B_h(p,s)$ are also euclidean balls $B = B_e(w,r)$, but
centers and radii vary: Namely, if $w=(a,b)$ is the euclidean center of $B$
and $r < b$ its euclidean radius, then the hyperbolic center of $B$ is
$p=(a, \sqrt{b^2-r^2})$, and the hyperbolic radius is
$s=\log\sqrt{\frac{b+r}{b-r}}$.

So as a topological space, our set $X$ is just $\{(x,y) \in \Bbb R^2: y>0\}$
with the usual topology, metrized by the hyperbolic distance. We define 
a Borel measure $\mu$ on $X = \Bbb H$ as follows: let $m_1$ be the restriction
to the upper half plane of the usual Gaussian probability
on $\Bbb R^2$, let $dm_2:= \chi _A dw$, where $dw$ stands
for planar Lebesgue measure and 
$A:= \{(x,y) \in \Bbb R^2: x > 1, x^{-1} > y>0\}$, and finally, set
$\mu := m_1 + m_2$.

First we show that 
$$M_\mu f(w) := \sup_{s >0}\frac1{\mu B_h(w,s)}\int_{B_h(w,s)}|f| d\mu$$
is not of weak type $(1,1)$.
By the usual approximation argument via convolutions, we
may use a Dirac delta instead of a function. So consider
$\delta_{(R+1/2,1)}$, where $R>>0$.
We will see that
$$\mu\left\{M_\mu \delta_{(R+1/2,1)}(w)
>\frac{(R-1)^{3/2}}3\right\}>\frac{1}{2R},$$ from which the result follows by
letting $R\to \infty$. Since 
$$
\mu\{(x,y) \in \Bbb R^2: R < x < R + 1, 
0 < y < x^{-1}\}
$$
$$
> m_2 \{(x,y) \in \Bbb R^2: R < x < R + 1, 
0 < y < x^{-1}\} =\log \left(1+\frac1R\right)>\frac1{2R},
$$
it is enough to prove that
$$\{(x,y) \in \Bbb R^2: R < x < R + 1, 
0 < y < x^{-1}\}\subset 
\left\{M_\mu \delta_{(R+1/2,1)}(w) >\frac{(R-1)^{3/2}}3\right\}.
$$
Fix $(x,y)$ such that $R < x < R + 1, 
0 < y < x^{-1}$, and select $r\in [1/2, 1)$ so that the
hyperbolic center of $B_e((x,1), r)$ is $(x,y)$. Since
$B_e((x,1), r)$ contains the point $(R+1/2,1)$,
$$M_\mu \delta_{(R+1/2,1)}((x,y)) \ge \frac1{\mu B_e((x,1),r)} \ge \frac1{\mu B_e((R,1),1)}.$$
To estimate $\mu B_e((R,1),1)$, note first that
$$m_1\ B_e((R,1),1) <
\int_{R-1}^\infty e^{-t^2/2} dt <
\int_{R-1}^\infty \frac{t}{R-1}e^{-t^2/2} dt =
 \frac{e^{-(R-1)^2/2}}{R-1}.$$
Suppose next that $(x-R)^2 + (y-1)^{2} <1$ and $0 < y < x^{-1}.$
Solving for $y$ in $(x-R)^2 + (y-1)^{2} <1$ we get
$$\frac{(x-R)^2}2<1-\sqrt{1 - (x-R)^2} < y < x^{-1} < \frac 1{R-1}.$$
So
$$\{(x,y) \in \Bbb R^2: (x-R)^2 + (y-1)^{2} <1,  0 < y < x^{-1}\}
$$
$$
\subset \left\{(x,y) \in \Bbb R^2: R - \sqrt{\frac 2{R-1}} < x < R + \sqrt{\frac 2{R-1}}, 
0 < y < (R-1)^{-1}\right\},$$
whence
$$m_2\ B_e((R,1),1) <
\int_0^{\frac1{R-1}}\int_{R - \sqrt{\frac2{R-1}}}^{R + \sqrt{\frac2{R-1}}} dx dy 
= \frac{2\sqrt 2}{(R-1)^{3/2}}.$$
Thus, by taking $R$ sufficiently large,
$$\mu\ B_e((R,1),1) <
 \frac{3}{(R-1)^{3/2}},$$
and it follows that $M_\mu$ is not of weak type $(1,1)$.

Next we prove that there is a constant $c> 0$ such that for
all $w\in \Bbb H$, $\mu B_h(w,s)\le c s$. 
We find constants $c_1, c_2 > 0$ with $m_1 B_h(w,s)\le c_1 s$ and 
 $m_2 B_h(w,s)\le c_2 s$ for
all $w\in \Bbb H$. From now
 on we shall adhere to the convention whereby a constant
 $c$ may change its value from one line to the next.

Recall that we use $B_e((a,b), r)= B_h((a,b^\prime),s)$ to denote the same 
ball $B\subset \Bbb H$, with respect to the euclidean metric in the first
case and the hyperbolic metric in the second.
Here
$$s = \frac12\log\left(1 + \frac{2r}{b-r}\right).$$
For small values of $b$, and therefore of $r$, $s$ controls
$r$, so $m_1$ and $m_2$ can simply be replaced by planar
Lebesgue measure to prove the polynomial growth. More precisely,
suppose $0<b\le 3$. Then $0<r<3$ and there exists a $c > 0$ such that
$s\ge c r$. So for $i=1,2$, 
$$m_i B_h((a,b^\prime),s) = m_i B_e((a,b), r)\le
c r\le c^\prime s.$$
Suppose next that $b >3$. Then 
$m_2\ B_e((a,b),r) = 0$ unless $r>b-1$. In this case
 we have
$$m_2 B_h((a,b^\prime),s) = m_2 B_e((a,b), r)\le
\log (1 + 2r)\le\log\left(1+\frac{2r}{b-r}\right) = 2s.$$

Finally, to prove the polynomial growth of $m_1$ on $\{b > 3\}$
we consider the following cases.
If $r$ and $b$ are comparable, say, $b/3 \le  r < b$, everything is
trivial, since 
$$m_1 B_h((a,b^\prime),s) = m_1 B_e((a,b), r)\le 1/2$$
and
$$s\ge\frac12\log 2 .$$
So suppose $0 < r <b/3$. Then
$$m_1 B_h((a,b^\prime),s) = m_1 B_e((a,b), r)
\le m_1 B_e((0,b), r)
$$
$$
\le\frac1{2\pi}\int_{-r}^{r}\int_{b-r}^\infty  e^{-\frac{x^2 + y^2}{2}} dy dx
\le\frac1{2\pi}\frac{e^{-\frac{(b-r)^2}2}}{b-r}\int_{-r}^{r}
e^{-\frac{x^2}{2}}  dx 
\le\frac{e^{-\frac{(b-r)^2}2}}{b-r}\min\{1,2r\}.
$$
Now if $0<r<1/2$, then
$$\frac{e^{-\frac{(b-r)^2}2}}{b-r}2r\le \frac{c2r}b\le c^\prime \frac12\log\left(1 + \frac{2r}{b}\right)\le c^\prime s,
$$
while if $1/2 \le r <b/3$, then
$$\frac{e^{-\frac{(b-r)^2}2}}{b-r}\le \frac{c^\prime}b\le c^{\prime\prime}
\frac12\log\left(1 + \frac{1}{b}\right)\le c^{\prime\prime} s.
$$
\qed

\noindent{\em Remark.} One might ask for which locally finite Borel
measures $\mu$  on $\Bbb H$ is the centered
maximal function $M_\mu$ of weak type $(1,1)$. Locally finite means that
for every $x\in \Bbb H$ there is an open neighborhood of $x$ with
finite measure. This implies that compact sets and balls have finite
measure. If 
$\mu$ has compact support (and hence it is finite), then it follows from
Besicovitch's covering theorem that
$M_\mu$ is of weak type $(1,1)$. But being finite, even in the presence
of polynomial growth, is not enough to ensure the weak type of $M_\mu$.
The example we present above can be easily modified so that $\mu (\Bbb H ) <
\infty$: Instead of using the set 
$\{(x,y) \in \Bbb R^2: x > 1, x^{-1} > y>0\}$ to define $m_2$, take
for instance
$A:= \{(x,y) \in \Bbb R^2: x > 0, e^{-x} > y>0\}$, and argue as before.

On the other hand, the centered 
maximal function associated to area
in the hyperbolic plane is of weak type $(1,1)$
(cf. \cite{Str}), even though area does not
satisfy any polynomial growth condition.
So there seems to
be no significant 
relationship between the weak type of the maximal operator and the polynomial
growth 
of the underlying measure.

To finish, we mention that while the doubling condition
on the measure is sufficient to ensure the weak type of the
maximal function, for $\Bbb H$ this is irrelevant: It follows from the remark
in page 67 of \cite{CW} together with  Example 3.5.2 of \cite{Luu} that the
hyperbolic plane admits no doubling measures (in particular, area is not
doubling).

\end{document}